\documentclass[12]{amsart}
\usepackage{amsfonts}
\usepackage{indentfirst,amsthm,bm,amsmath}
\date{}
\newtheorem{theorem}{Theorem}[section]
\newtheorem{corollary}{Corollary}[section]
\newtheorem{lemma}{Lemma}[section]

\begin{document}
\title[Uniqueness theorems for meromorphic mappings]{Uniqueness theorems for meromorphic mappings
sharing hyperplanes in general position}
\thanks{This work was partially supported by the NSFC (No. 10771121), the NSF of Jiangxi (No. 2008GQS0075) and the YFED of
Jiangxi (No. GJJ10050) of China.}
\author{Ting-Bin Cao \and Hong-Xun Yi}
\address[Ting-Bin Cao]{Department of Mathematics, Nanchang University, Nanchang,
Jiangxi 330031, China} \email{tbcao@ncu.edu.cn, ctb97@163.com (the
corresponding author)}
\address[Hong-Xun Yi]{Department of Mathematics, Shandong University, Jinan,
Shandong 250100, China} \email{hxyi@sdu.edu.cn}
\date{}
\keywords{meromorphic mapping, uniqueness theorem, Nevanlinna
theory, hyperplane} \subjclass[2000]{32H30, 32A22, 30D35}
\begin{abstract}The purpose of this article is to study the uniqueness problem for
meromorphic mappings from $\mathbb{C}^{n}$ into the complex
projective space $\mathbb{P}^{N}(\mathbb{C}).$ By making using of
the method of dealing with multiple values due to L. Yang and the
technique of Dethloff-Quang-Tan respectively, we obtain two general
uniqueness theorems which improve and extend some known results of
meromorphic mappings sharing hyperplanes in general
position.\end{abstract} \maketitle

\section{Introduction and main results}
For a nonconstant meromorphic function $f$ on $\mathbb{C}^{n}$ and
$a\in\mathbb{P}^{1}(\mathbb{C}),$ we denote by $\nu_{f-a}$ the map
from $\mathbb{C}$ into $\mathbb{Z}$ whose value $\nu_{f-a}(z)$ is
the multiplicity of the zero of $f-a$ at $z.$\par

In 1926, R. Nevanlinna \cite{nevanlinna} proved the well-known
five-value theorem that for two nonconstant meromorphic functions
$f$ and $g$ on the complex plane $\mathbb{C},$ if they have the same
inverse images (ignoring multiplicities) for five distinct values in
$\mathbb{P}^{1}(\mathbb{C}),$ then $f(z)\equiv g(z).$ We know that
the number five of distinct values in Nevanlinna's five-value
theorem cannot be reduced to four. For example, $f(z)=e^{z}$ and
$g(z)=e^{-z}$ share four values $0, 1, -1, \infty$ (ignoring
multiplicities), but $f(z)\not\equiv g(z).$  There have been several
improvements of Nevanlinna's five-value theorem. H. X. Yi
(\cite{yi-yang}, Theorem 3.15) adopted the method of dealing with
multiple values due to L. Yang \cite{yang} and obtained a uniqueness
theorem of meromorphic functions of one variable. Later, Hu, Li and
Yang extended this result to meromorphic functions in several
variables (see Theorem 3.9 in \cite{hu-li-yang}).\par

\begin{theorem}(\cite{hu-li-yang}, Theorem 3.9) \label{T-B} Let $f$ and
$g$ be two nonconstant meromorphic functions on $\mathbb{C}^{n},$
let $a_{j}$ $(j=1,2,\ldots,q)$ be $q$ distinct complex elements in
$\mathbb{P}^{1}(\mathbb{C})$ and take
$m_{j}\in\mathbb{Z}^{+}\cup\{\infty\}$ $(j=1,2,\ldots,q)$ satisfying
$m_{1}\geq m_{2}\geq\cdots\geq m_{q}$ and $\nu_{f-a_{j}, \leq
m_{j}}^{1}=\nu_{g-a_{j}, \leq m_{j}}^{1}$ $(j=1,2,\ldots,q).$ If
$\sum_{j=3}^{q}\frac{m_{j}}{m_{j}+1}>2,$ then $f(z)\equiv
g(z).$\end{theorem}

Over the last few decades, there have been several generalizations
of Nevanlinna's five-value theorem to the case of meromorphic
mappings from $\mathbb{C}^{n}$ into the complex projective space
$\mathbb{P}^{N}(\mathbb{C}).$ Some of the first results concerning
this research are due to Fujimoto \cite{fujimoto-1, fujimoto-2}.
\par

For a meromorphic mapping $f$ from $\mathbb{C}^{n}$ into
$\mathbb{P}^{N}(\mathbb{C})$ and a hyperplane $H$ in
$\mathbb{P}^{N}(\mathbb{C}),$ we denote by $\nu_{(f, H)}$ the map
from $\mathbb{C}$ into $\mathbb{Z}$ whose value $\nu_{(f, H)}(z)$
$(z\in \mathbb{C}^{n})$ is the intersection multiplicity of the
images of $f$ and $H$ at $f(z).$ Let $H_{1},$ $H_{2},$ $\ldots,$
$H_{q}$ be $q$ hyperplanes in general position such that $dim
f^{-1}(H_{i}\cap H_{j})\leq n-2$ for $i\neq j.$ Take $d, m_{j}$ be
positive integers or $\infty.$ Consider the set $\mathcal{F}_{\leq
m_{j}}(f, \{H_{j}\}_{j=1}^{q}, d)$ of all linearly non-degenerate
meromorphic mappings $g:
\mathbb{C}^{n}\rightarrow\mathbb{P}^{N}(\mathbb{C})$ satisfying the
conditions:\par

(a)$\nu_{(f, H_{j}), \leq m_{j}}^{d}=\nu_{(g, H_{j}), \leq
m_{j}}^{d},$\par

(b) $f(z)=g(z)$ on $\cup_{j=1}^{q}\{z\in\mathbb{C}^{n}| 0<\nu_{(f,
H_{j})}\leq m_{j}\}.$\\
For brevity we will omit the notation $\leq m_{j}$ if
$m_{j}=\infty.$\par

Fujimoto \cite{fujimoto-1, fujimoto-2} proved that if $q\geq 3N+2,$
then $\sharp\mathcal{F}(f, \{H_{j}\}_{j=1}^{q}, \infty)=1.$ In 1983,
Smiley \cite{smiley} obtained an improvement with truncated number
one that if $q\geq 3N+2,$ then $\sharp\mathcal{F}(f,
\{H_{j}\}_{j=1}^{q}, 1)=1.$ Later, Thai and Quang \cite{thai-quang}
considered a smaller number $q$ and proved that if $N\geq 2$ and
$q\geq 3N+1,$ then $\sharp\mathcal{F}(f, \{H_{j}\}_{j=1}^{q}, 1)=1.$
In \cite{dethloff-tan}, Dethloff and Tan considered the case $q\geq
2N+3$ and obtained that if $q\geq 2N+3,$ then $\sharp\mathcal{F}(f,
\{H_{j}\}_{j=1}^{q}, N)=1.$ Recently, Chen and Yan \cite{chen-yan-0}
improved the above results and obtained that if $q\geq 2N+3,$ then
$\sharp\mathcal{F}(f, \{H_{j}\}_{j=1}^{q}, 1)=1.$\par

Considering multiple values, there are some theorems. The following
result is Theorem 0.1 in \cite{aihara} for the special case of
linearly non-degenerate meromorphic mappings.\par

\begin{theorem} \cite{aihara}\label{T-C} Let $f$ and $g$ be two linearly non-degenerate meromorphic mappings
of $\mathbb{C}^{n}$ into $\mathbb{P}^{N}(\mathbb{C}),$ and let
$H_{1},$ $H_{2},$ $\ldots,$ $H_{q}$ be $q$ hyperplanes in general
position such that  $dim f^{-1}(H_{i}\cap H_{j})\leq n-2$ for $i\neq
j.$ Take $m_{j} (j=1,2,\ldots,q)$ be positive integers or $\infty$
such that $m_{1}\geq m_{2}\geq\cdots\geq m_{q}\geq N,$
$$\nu^{1}_{(f, H_{j}), \leq m_{j}}=\nu^{1}_{(g, H_{j}), \leq
m_{j}}\,\,(j=1,2,\ldots,q),$$ and $f(z)=g(z)$ on
$\cup_{j=1}^{q}\{z\in\mathbb{C}^{n}| 0<\nu_{(f, H_{j})}\leq
m_{j}\}.$ If
\begin{equation}\label{E-0}\sum_{j=3}^{q}\frac{m_{j}}{m_{j}+1}>\frac{Nq-q+N+1}{N}
+\frac{m_{1}}{m_{1}+1}-\frac{m_{2}}{m_{2}+1},\end{equation}
then $f(z)\equiv g(z).$\end{theorem}

Remark that condition \eqref{E-0} implies $q\geq 3N+2.$ In 2006,
Dethloff and Tan \cite{dethloff-tan-1} obtained the following result
for a smaller $q.$\par

\begin{theorem} \cite{dethloff-tan-1}\label{T-A} If $N\geq 7,$ then $\sharp\mathcal{F}_{\leq m}(f, \{H_{j}\}_{j=1}^{3N-1},
1)=1,$ where $m\geq 6N+32+\frac{192}{N-6}.$ If $4\geq N\geq 6,$ then
$\sharp\mathcal{F}_{\leq m}(f, \{H_{j}\}_{j=1}^{3N}, 2)=1,$ where
$m\geq 6N+17+\frac{51}{N-3}.$
\end{theorem}

The main purpose of this paper is to consider the multiple values
and uniqueness problem of meromorphic mappings. By making use of the
method of dealing with multiple values due to L. Yang \cite{yang},
we obtain the first main result below.\par

\begin{theorem}\label{T-0} Let $f$ and $g$ be two linearly non-degenerate meromorphic mappings
of $\mathbb{C}^{n}$ into $\mathbb{P}^{N}(\mathbb{C}),$ and let
$H_{1},$ $H_{2},$ $\ldots,$ $H_{q}$ be $q$ hyperplanes in general
position such that  $dim f^{-1}(H_{i}\cap H_{j})\leq n-2$ for $i\neq
j.$ Take $m_{j} (j=1,2,\ldots,q)$ be positive integers or $\infty$
such that $m_{1}\geq m_{2}\geq\cdots\geq m_{q}\geq N,$
$$\nu^{1}_{(f, H_{j}), \leq m_{j}}=\nu^{1}_{(g, H_{j}), \leq
m_{j}}\,\,(j=1,2,\ldots,q),$$ and $f(z)=g(z)$ on
$\cup_{j=1}^{q}\{z\in\mathbb{C}^{n}| 0<\nu_{(f, H_{j})}\leq
m_{j}\}.$ If
\begin{equation}\label{E-1}\sum_{j=3}^{q}\frac{m_{j}}{m_{j}+1}>\frac{Nq-q+N+1}{N}-
(N-1)\left(\frac{1}{m_{1}+1}+\frac{1}{m_{2}+1}\right),\end{equation}
then $f(z)\equiv g(z).$\end{theorem}

Obviously, Theorem \ref{T-0} is an improvement of Theorem \ref{T-C}.
For the case $N=1,$ the condition \eqref{E-1} reduces to
$\sum_{j=3}^{q}\frac{m_{j}}{m_{j}+1}>2.$ Thus Theorem \ref{T-0} is
an extension of Theorem \ref{T-B}. Furthermore, we immediately get
the following corollaries from which we see that Theorem \ref{T-0}
is an improvement of Smiley's $3N+2$ hyperplanes uniqueness theorem
\cite{smiley}.\par

\begin{corollary}\label{C-1} If $q\geq 3N+2,$ then $\sharp\mathcal{F}_{\leq m}(f,
\{H_{j}\}_{j=1}^{q}, 1)=1,$ where $m$ $>$
$N-1+\frac{N(N+1)}{q-3N-1}.$
\end{corollary}

\begin{corollary}\label{C-2} $\sharp\mathcal{F}_{\leq m}(f,
\{H_{j}\}_{j=1}^{3N+2}, 1)=1,$ where $m>N^{2}+2N-1.$
\end{corollary}

Considering a smaller $q$ than $3N+2,$ we have another main theorem
by using the technique shown in \cite{chen-yan-0,
dethloff-quang-tan, dethloff-tan}.\par

\begin{theorem}\label{T-1} Let $f$ and $g$ be two linearly non-degenerate meromorphic mappings
of $\mathbb{C}^{n}$ into $\mathbb{P}^{N}(\mathbb{C}),$ and let
$H_{1},$ $H_{2},$ $\ldots,$ $H_{q}$ be $q$ $(q\geq 2N)$ hyperplanes
in general position such that  $dim f^{-1}(H_{i}\cap H_{j})\leq n-2$
for $i\neq j.$ Take $m_{j} (j=1,2,\ldots,q)$ be positive integers or
$\infty$ such that $m_{1}\geq m_{2}\geq\cdots\geq m_{q}\geq N,$
$$\nu^{1}_{(f, H_{j}), \leq m_{j}}=\nu^{1}_{(g, H_{j}), \leq m_{j}}\,\,(j=1,2,\ldots,q),$$ and $f(z)=g(z)$ on
$\cup_{j=1}^{q}\{z\in\mathbb{C}^{n}| 0<\nu_{(f, H_{j})}\leq
m_{j}\}.$ If
\begin{eqnarray}\label{E-2} \sum_{j=3}^{q}\frac{m_{j}}{m_{j}+1}&>&\frac{Nq-q+N+1}{N}-\frac{4N-4}{q+2N-2}
+\left(\frac{1}{m_{1}+1}+\frac{1}{m_{2}+1}\right),\end{eqnarray}
then $f(z)\equiv g(z).$\end{theorem}

Thus we obtain immediately the following corollaries which improve
the above-mentioned uniqueness theorems for meromorphic mappings
sharing hyperplanes in general position \cite{fujimoto-1, smiley,
thai-quang, dethloff-tan, dethloff-tan-1, chen-yan-0}.\par

\begin{corollary}\label{C-3}
If $q\geq 2N+3,$ then $\sharp\mathcal{F}_{\leq m}(f,
\{H_{j}\}_{j=1}^{q}, 1)=1,$ where
$$m>\frac{(N-1)q^{2}+(2N^{2}-N+3)q+2N^{2}-2}{(q+N-1)(q-2N-2)}.$$
\end{corollary}

\begin{corollary}\label{C-4}$\sharp\mathcal{F}_{\leq m}(f,
\{H_{j}\}_{j=1}^{2N+3}, 1)=1,$ where
$m>\frac{8N^{3}+14N^{2}-2}{3N+2}.$
\end{corollary}

\begin{corollary}\label{C-5}If $N\geq 2,$ then $\sharp\mathcal{F}_{\leq m}(f,
\{H_{j}\}_{j=1}^{3N+1}, 1)=1,$ where
$m>\frac{15N^{2}-2N+3}{4(N-1)}.$\end{corollary}

\begin{corollary} \label{C-6} If $N\geq 3,$ then
$\sharp\mathcal{F}_{\leq m}(f, \{H_{j}\}_{j=1}^{3N}, 1)=1,$ where
$m>\frac{15N^{3}-10N^{2}+9N-2}{(4N-1)(N-2)}.$ \end{corollary}

\begin{corollary} \label{C-7} If $N\geq 4,$ then
$\sharp\mathcal{F}_{\leq m}(f, \{H_{j}\}_{j=1}^{3N-1}, 1)=1,$ where
$m>\frac{15N^{3}-18N^{2}+17N-6}{(4N-2)(N-3)}.$
\end{corollary}

The last corollary is a supplement of Corollary \ref{C-2}.\par

\begin{corollary}\label{C-8}$\sharp\mathcal{F}_{\leq m}(f,
\{H_{j}\}_{j=1}^{3N+2}, 1)=1,$ where $m>\frac{15N^{2}+6N-1}{4N+1}.$
\end{corollary}

For the case $N=1,$ the condition \eqref{E-2} reduces to
$\sum_{j=3}^{q}\frac{m_{j}}{m_{j}+1}>2+\left(\frac{1}{m_{1}+1}+\frac{1}{m_{2}+1}\right).$
Thus compared with the conditions of Theorems \ref{T-B} and
\ref{T-0}, there maybe exist a better lower estimate than condition
\eqref{E-2} in Theorem \ref{T-1}.\par

\section{Preliminaries}
We set $\|z\|=(\sum_{j=1}^{n}|z_{j}|^{2})^{\frac{1}{2}}$ for
$z=(z_{1}, \ldots, z_{n})\in\mathbb{C}^{n}.$ For $r>0,$ define
$B(r)=\{z\in\mathbb{C}^{n}: \|z\|<r\},$ $S(r)=\{z\in\mathbb{C}^{n}:
\|z\|=r\},$
$d^{c}=(4\pi\sqrt{-1})^{-1}(\partial-\overline{\partial}),$
$\upsilon=(dd^{c}\|z\|^{2})^{n-1}$ and
$\sigma=d^{c}\log\|z\|^{2}\wedge(dd^{c}\|z\|^{2})^{n-1}.$\par

Let $h$ be a nonzero entire function on $\mathbb{C}^{n}.$ For $a\in
\mathbb{C}^{n},$ we can write $h$ as
$h(z)=\sum_{m=0}^{\infty}P_{m}(z-a),$ where the term $P_{m}(z)$ is
either identically zero or a homogeneous polynomial of degree $m.$
The number $\nu_{h}(a):=\min\{m: P_{m}\neq 0\}$ is said to be the
zero-multiplicity of $h$ at $a.$ Set
$supp\nu_{h}:=\overline{\{z\in\mathbb{C}^{n}: \nu_{h}(z)\neq
0\}}.$\par

Let $\varphi$ be a nonzero meromorphic function on $\mathbb{C}^{n}.$
For each $a\in\mathbb{C}^{n},$ we choose nonzero holomorphic
functions $\varphi_{0}$ and $\varphi_{1}$ on a neighborhood $U$ of
$a$ such that  $\varphi=\frac{\varphi_{0}}{\varphi_{1}}$ on $U$ and
$dim (\varphi_{0}^{-1}\cap \varphi_{1}^{-1}(0))\leq n-2,$ and we
define $\nu_{\varphi}:=\nu_{\varphi_{0}},$
$\nu_{\varphi}^{\infty}:=\nu_{\varphi_{1}},$ which are independent
of choices of $\varphi_{0}$ and $\varphi_{1}.$ \par

Let $f$ be a nonconstant meromorphic mapping of $C^{n}$ into
$\mathbb{P}^{N}(\mathbb{C}).$ We can choose holomorphic functions
$f_{0},$ $f_{1},$ $\ldots,$ $f_{N}$ on $\mathbb{C}^{n}$ such that
$I_{f}:=\{z\in\mathbb{C}^{n}: f_{0}(z)=\cdots=f_{N}(z)=0\}$ is of
dimension at most $n-2$ and $f=(f_{0}: \cdots: f_{N}).$ Usually,
$(f_{0}: \cdots: f_{N})$ is a reduces representation of $f.$ The
characteristic function of $f$ is defined by
\begin{equation*}
    T(r, f)=\int_{S(r)}\log\|f\|\sigma-\int_{S(1)}\log\|f\|\sigma\quad (r>1),
\end{equation*}
where $\|f\|=(\sum_{j=0}^{N}|f_{j}|^{2})^{\frac{1}{2}}.$ Note that
$T(r,f)$ is independent of the choice of the reduced representation
of $f$.\par

For a hyperplane $H=\{(x_{0}: \cdots:
x_{N})\in\mathbb{P}^{N}(\mathbb{C}):
a_{0}x_{0}+\cdots+a_{N}x_{N}=0\},$ we say that $(f, H)$ is free if
$(f, H)=\sum_{j=0}^{N}a_{j}f_{j}\not\equiv 0.$ Under the assumption
that $(f, H)$ is free, $(f, H)$ is a nonzero holomorphic function
and the proximity function of $f$ and $H$ is defined by

$$m_{f, H}(r)=\int_{S(r)}\log\frac{\|f\|\|H\|}{|(f, H)|}\sigma-\int_{S(1)}
\log\frac{\|f\|\|H\|}{|(f, H)|}\sigma, \,\, r>1,$$ where
$\|H\|=(\sum_{j=0}^{N}|a_{j}|^{2})^{\frac{1}{2}}.$ The proximity
function of a meromorphic function $\varphi$ on $\mathbb{C}^{n}$ is
defined by  $m(r, \varphi)=\int_{S(r)}\log^{+}|\varphi|\sigma,$
where $\log^{+}x=\max\{\log x, 0\}$ for $x\geq 0.$\par

Let $k, M$ be positive integers or $+\infty.$ For a divisor $\nu$ on
$\mathbb{C}^{n}.$ We define the counting functions of $\nu$ as
follows. Set
\begin{equation*}
\nu^{M}(z)=\min\{\nu(z), M\},\quad \nu_{\leq k}^{M}(z)=
\left\{\begin{array}{ll}
 0, & \hbox{if}\quad \nu(z)>k; \\
 \nu^{M}(z), & \hbox{if}\quad \nu(z)\leq k,
\end{array}
\right.
\end{equation*} \begin{equation*} \nu_{\geq k}^{M}(z)= \left\{\begin{array}{ll}
 0, & \hbox{if}\quad \nu(z)< k; \\
 \nu^{M}(z), & \hbox{if}\quad \nu(z)\geq k.
\end{array}
\right.
\end{equation*} and
\begin{eqnarray*} n(t)=\left\{
                                    \begin{array}{ll}
                                      \int_{supp\nu\cap B(t)}\nu(z)\upsilon, & \hbox{if $n\geq 2;$} \\
                                      \sum_{|z|\leq t}\nu(z), & \hbox{if $n=1.$}
                                    \end{array}
                                  \right.
\end{eqnarray*}Similarly, we define $n^{M}(t),$ $n_{\geq k}^{M}(t)$ and $n_{\leq k}^{M}(t).$ We define
\begin{equation*}
N(r, \nu)=\int_{1}^{r}\frac{n(t)}{t^{2n-1}}dt\quad
(r>1).\end{equation*} Similarly, we define $N(r, \nu^{M}),$ $N(r,
\nu^{M}_{\leq k})$ and $N(r, \nu^{M}_{\geq k})$ and denote them by
$N^{M}(r, \nu),$ $N^{M}_{\leq k}(r, \nu)$ and $N^{M}_{\geq k}(r,
\nu),$ respectively.\par

For a meromorphic function $\varphi$ on $\mathbb{C}^{n},$ we denote
by $$N_{\varphi}(r)=N(r, \nu_{\varphi}),\quad
N^{M}_{\varphi}(r)=N^{M}(r, \nu_{\varphi}),$$
$$ N^{M}_{\varphi, \leq
k}(r)=N^{M}_{\leq k}(r, \nu_{\varphi}),\quad N^{M}_{\varphi, \geq
k}(r)=N^{M}_{\geq k}(r, \nu_{\varphi}).$$ For brevity we will omit
the superscript $M$ if $M=\infty.$\par
\par

For a closed subset $A$ of a purely $(n-1)$-dimensional analytic
subset of $\mathbb{C}^{n},$ we define
\begin{eqnarray*} n^{1}_{A}(t)=\left\{
                                    \begin{array}{ll}
                                      \int_{A\cap B(t)}\upsilon, & \hbox{if $n\geq 2;$} \\
                                      \sharp(A\cap B(t)), & \hbox{if $n=1,$}
                                    \end{array}
                                  \right.
\end{eqnarray*}and \begin{equation*}
N^{1}(r, A)=\int_{1}^{r}\frac{n^{1}_{A}(t)}{t^{2n-1}}dt\quad
(r>1).\end{equation*}

The First Main Theorem is that
\begin{equation*}
 T(r, f)=N_{(f, H)}(r)+m_{f, H}(r)+O(1).
\end{equation*}

As usual, by the notation $"\| P"$ we mean the assertion $P$ holds
for all $r>1$ excluding a set of finite Lebesgue measure.\par

\begin{theorem}(Second Main Theorem) Let $f$ be a linearly non-degenerate meromorphic mapping
from $\mathbb{C}^{n}$ into $\mathbb{P}^{N}(\mathbb{C})$ and $H_{1},$
$\ldots,$ $H_{q}$ be $q$ hyperplanes in general position. Then
\begin{equation*}\|(q-N-1)T(r, f)\leq\sum_{j=1}^{q}N^{N}_{(f, H_{j})}(r)+o(T(r,f)).\end{equation*}
 \end{theorem}

For two distinct hyperplanes $H_{1}$ and $H_{2},$ we have

\begin{lemma}\cite{ru}\label{L-1} $T(r, \frac{(f, H_{1})}{(f, H_{2})})\leq T(r, f)+O(1).$\end{lemma}

The following lemma is proved by using the method of dealing with
multiple values due to L. Yang \cite{yang}, see also Lemma 4.7 in
\cite{thai-quang}.\par

\begin{lemma} \label{L-2} Let $f$ be a linearly non-degenerate meromorphic  mapping
from $\mathbb{C}^{n}$ into $\mathbb{P}^{N}(\mathbb{C}),$ $H$ be a
hyperplane in general position, and $k(\geq N)$ be a positive
integer. Then
\begin{eqnarray*}
N^{N}_{(f, H)}(r)&\leq&N(1-\frac{N}{k+1})N^{1}_{(f, H),\leq k}(r)
+\frac{N}{k+1}N_{(f, H)}(r);
\end{eqnarray*}
and
\begin{eqnarray*}
N^{N}_{(f, H)}(r)&\leq&N(1-\frac{N}{k+1})N^{1}_{(f, H),\leq k}(r)
+\frac{N}{k+1}T(r,f)+o\left(T(r,f)\right).
\end{eqnarray*}
 \end{lemma}

\begin{proof} From
\begin{equation*}
N^{N}_{(f, H)}(r)=N^{N}_{(f, H),\leq k}(r)+N^{N}_{(f, H),\geq
k+1}(r)
\end{equation*}
and
\begin{equation*}
N^{N}_{(f, H),\geq k+1}(r)\leq\frac{N}{k+1}N_{(f, H),\geq
k+1}(r)\leq\frac{N}{k+1}\left(N_{(f, H)}(r)-N^{N}_{(f,H),\leq
k}(r)\right),
\end{equation*}
we deduce that
\begin{eqnarray*}
N^{N}_{(f, H)}(r)&\leq&(1-\frac{N}{k+1})N^{N}_{(f, H),\leq k}(r)
+\frac{N}{k+1}N_{(f, H)}(r)\\&\leq&N(1-\frac{N}{k+1})N^{1}_{(f,
H),\leq k}(r) +\frac{N}{k+1}N_{(f, H)}(r).
\end{eqnarray*}This completes the proof of the first inequality of
the lemma. The second inequality of the lemma follows immediately
because of $N_{(f, H)}(r)\leq T(r,f)+o\left(T(r,f)\right).$
\end{proof}\par\vskip 6pt

\section{Proof of Theorem \ref{T-0}}
 Suppose that
$f(z)\not\equiv g(z).$ By the Second Main Theorem, we have
\begin{eqnarray*}
\left\|(q-N-1)T(r)\right.&\leq&\sum_{j=1}^{q}\left(N^{N}_{(f,
H_{j})}(r)+N^{N}_{(g, H_{j})}(r)\right) +o\left(T(r)\right),
\end{eqnarray*}where $T(r)=T(r,f)+T(r, g).$
By Lemma \ref{L-2} we have
\begin{eqnarray*}
N^{N}_{(f,H_{j})}(r)+N^{N}_{(g,H_{j})}(r)&\leq&N\left(1-\frac{N}{m_{j}+1}\right)\left(N^{1}_{(f,H_{j}),\leq
m_{j}}(r)+N^{1}_{(g,H_{j}),\leq
m_{j}}(r)\right)\\&&+\frac{N}{m_{j}+1}T(r)+o\left(T(r)\right).
\end{eqnarray*}
Therefore,  we have
\begin{eqnarray*}
\left\|(q-N-1)T(r)\right. &\leq&
N\sum_{j=1}^{q}\left(1-\frac{N}{m_{j}+1}\right)\left(N^{1}_{(f,
H_{j}), \leq m_{j}}(r)+N^{1}_{(g, H_{j}), \leq
m_{j}}(r)\right)\\&&+N\sum_{j=1}^{q}
\frac{1}{m_{j}+1}T(r)+o\left(T(r)\right).
\end{eqnarray*}\par

Noting that $m_{1}\geq m_{2}\geq\cdots\geq m_{q}\geq N,$ we have
\begin{eqnarray*}&&\sum_{j=1}^{q}\left(1-\frac{N}{m_{j}+1}\right)N^{1}_{(f,
H_{j}),\leq m_{j}}(r)\\&=&\left(1-\frac{N}{m_{1}+1}\right)N^{1}_{(f,
H_{1}), \leq
m_{1}}(r)+\sum_{j=2}^{q}\left(1-\frac{N}{m_{j}+1}\right)N^{1}_{(f,
H_{j}), \leq
m_{j}}(r)\\&\leq&\left(1-\frac{N}{m_{1}+1}\right)N^{1}_{(f, H_{1}),
\leq
m_{1}}(r)+\sum_{j=2}^{q}\left(1-\frac{N}{m_{2}+1}\right)N^{1}_{(f,
H_{j}), \leq
m_{j}}(r)\\&=&\left(\frac{N}{m_{2}+1}-\frac{N}{m_{1}+1}\right)N^{1}_{(f,
H_{1}), \leq
m_{1}}(r)+\sum_{j=1}^{q}\left(1-\frac{N}{m_{2}+1}\right)N^{1}_{(f,
H_{j}), \leq
m_{j}}(r)\\&\leq&\left(\frac{N}{m_{2}+1}-\frac{N}{m_{1}+1}\right)T(r,f)
+\sum_{j=1}^{q}\left(1-\frac{N}{m_{2}+1}\right)N^{1}_{(f, H_{j}),
\leq m_{j}}(r)+O(1).\end{eqnarray*} The same is true for $g$ instead
of $f.$ Hence, we can deduce that
\begin{eqnarray*}
&&\left\|(q-N-1)T(r)\right.\\&\leq&
N^{2}\left(\frac{1}{m_{2}+1}-\frac{1}{m_{1}+1}\right)T(r)+N\sum_{j=1}^{q}
\frac{1}{m_{j}+1}T(r)\\&&+N\left(1-\frac{N}{m_{2}+1}\right)\sum_{j=1}^{q}\left(N^{1}_{(f,
H_{j}), \leq m_{j}}(r)+N^{1}_{(g, H_{j}), \leq
m_{j}}(r)\right)+o\left(T(r)\right).
\end{eqnarray*}

Since $\nu_{(f, H_{j}), \leq m_{j}}^{1}=\nu_{(g, H_{j}), \leq
m_{j}}^{1}$ and $f(z)=g(z)$ on $\cup_{j=1}^{q}\{z\in\mathbb{C}^{n}|
0<\nu_{(f, H_{j})}\leq m_{j}\},$ we have
$$\sum_{j=1}^{q}\left(N^{1}_{(f, H_{j}),\leq
m_{j}}(r)+N^{1}_{(g, H_{j}),\leq m_{j}}(r)\right)\leq
2N_{f-g}(r)\leq 2T(r)+O(1).$$ Hence, we can deduce that
\begin{eqnarray*}
\left\|(q-N-1)T(r)\right.&\leq& \left(N\sum_{j=1}^{q}
\frac{1}{m_{j}+1}+2N-N^{2}\left(\frac{1}{m_{1}+1}+\frac{1}{m_{2}+1}\right)\right)T(r)\\&&+o\left(T(r)\right).
\end{eqnarray*}
Noting that $q=\sum_{j=1}^{q}\frac{m_{j}+1}{m_{j}+1},$ we can obtain
from the above inequality  that
\begin{eqnarray*}\left\|\left(\sum_{j=3}^{q}\frac{m_{j}}{m_{j}+1}-\frac{Nq-q+N+1}{N}+
(N-1)\left(\frac{1}{m_{1}+1}+\frac{1}{m_{2}+1}\right)\right)T(r)\right.
\leq o(T(r)).\end{eqnarray*} This is a contradiction.

\section{Proof of Theorem \ref{T-1}}
 For brevity we denote $T(r, f)+T(r, g)$ by $T(r).$ Suppose that
$f(z)\not\equiv g(z).$ Then by changing indices, if necessary, we
may assume that
$$\underbrace{\frac{(f, H_{1})}{(g, H_{1})}\equiv \frac{(f, H_{2})}{(g, H_{2})}\equiv\cdots
\frac{(f, H_{k_1})}{(g, H_{k_1})}}_{\hbox{group 1}}
\not\equiv\underbrace{\frac{(f, H_{k_1+1})}{(g,
H_{k_1+1})}\equiv\cdots\equiv\frac{(f, H_{k_2})}{(g,
H_{k_2})}}_{\hbox{group
2}}$$$$\not\equiv\cdots\not\equiv\underbrace{\frac{(f,
H_{k_{s-1}+1})}{(g, H_{k_{s-1}+1})}\equiv\cdots\equiv\frac{(f,
H_{k_s})}{(g, H_{k_s})}}_{\hbox{group s}},$$ where $k_{s}=q.$ Then
the number of elements of every group is at most $N$ because
$f(z)\not\equiv g(z).$\par

The map $\sigma:$ $\{1, 2, \ldots, q\}\rightarrow\{1, 2, \ldots,
q\}$ is defined by \begin{equation*}\sigma(i)=\left\{
                     \begin{array}{ll}
                      i+N, & \hbox{if $i+N\leq q;$} \\
                       i+N-q, & \hbox{if $i+N>q.$}
                     \end{array}
                   \right.\end{equation*}
Obviously, $\sigma$ is bijective. Since $q\geq 2N,$
$|\sigma(i)-i|\geq N.$ Thus, $\frac{(f, H_{i})}{(g, H_{i})}$ and
$\frac{(f, H_{\sigma(i)})}{(g, H_{\sigma(i)})}$ belongs to distinct
group, and so $\frac{(f, H_{i})}{(g, H_{i})}\not\equiv \frac{(f,
H_{\sigma(i)})}{(g, H_{\sigma(i)})}.$ \par

We here use the technique in \cite{dethloff-quang-tan}, see also
\cite{chen-yan-0, dethloff-tan}. Set $P_{i}:=(f, H_{i})(g,
H_{\sigma(i)})-(f, H_{\sigma(i)})(g, H_{i})\not\equiv 0,$ where
$1\leq i\leq q.$ By the assumption and the definition of function
$P_{i}$ we get that for $k\in\{i, \sigma(i)\}$ and any
$z_{0}\in\{z\in| 1\leq \nu_{(f, H_{k})}\leq m_{k}\}$ $(=\{z\in|
1\leq \nu_{(g, H_{k})}\leq m_{k}\}),$ $z_{0}$ is a zero of $P_{i}$
with
$$\nu_{P_{i}}(z_{0})\geq \min\{\nu_{(f, H_{k})}(z_{0}), \nu_{(g,
H_{k})}(z_{0})\}$$ outside an analytic set of codimension $\geq 2.$
On the other hand, since $\nu^{1}_{(f, H_{k}), \leq
m_{k}}=\nu^{1}_{(g, H_{k}), \leq m_{k}}$ we have
\begin{equation*}\min\{\nu_{(f, H_{k})}(z_{0}),\nu_{(g,
H_{k})}(z_{0})\}\geq \nu^{N}_{(f, H_{k}), \leq
m_{k}}(z_{0})+\nu^{N}_{(g, H_{k}), \leq m_{k}}(z_{0})-N\nu^{1}_{(f,
H_{k}), \leq m_{k}}(z_{0}).
\end{equation*}
We also get that for any $j\in\{1, 2, \ldots, q\}\setminus\{i,
\sigma(i)\},$ any zero of $(f, H_{j})$ is also a zero of $P_{i}$
outside an analytic set of codimension $\geq 2.$ Thus we have
\begin{eqnarray*}\nu_{P_{i}}&\geq&\nu^{N}_{(f, H_{i}), \leq
m_{i}}+\nu^{N}_{(f, H_{\sigma(i)}), \leq m_{\sigma(i)}}+\nu^{N}_{(g,
H_{i}), \leq m_{i}}+\nu^{N}_{(g, H_{\sigma(i)}), \leq
m_{\sigma(i)}}\\&&-N\nu^{1}_{(f, H_{i}),\leq m_{i}}-N\nu^{1}_{(f,
H_{\sigma(i)}),\leq m_{\sigma(i)}}+\sum_{j=1, j\neq i,
\sigma(i)}^{q}\nu^{1}_{(f, H_{j}), \leq m_{j}}
\end{eqnarray*}outside an analytic set of codimension $\geq 2.$
Hence, for all $i\in\{1, 2, \ldots, q\}$ we have
\begin{eqnarray*}N_{P_{i}}&\geq&N^{N}_{(f, H_{i}),
\leq m_{i}}(r)+N^{N}_{(f, H_{\sigma(i)}), \leq
m_{\sigma(i)}}(r)+N^{N}_{(g, H_{i}), \leq m_{i}}(r)+N^{N}_{(g,
H_{\sigma(i)}), \leq m_{\sigma(i)}}(r)\\&&-NN^{1}_{(f, H_{i}),\leq
m_{i}}(r)-NN^{1}_{(f, H_{\sigma(i)}),\leq
m_{\sigma(i)}}(r)+\sum_{j=1, j\neq i, \sigma(i)}^{q}N^{1}_{(f,
H_{j}), \leq m_{j}}(r).
\end{eqnarray*}\par

On the other hand, by Jensen's formula we have \begin{eqnarray*}
N_{P_{i}}(r)&=&\int_{S(r)}\log|P_{i}|\sigma+O(1)\\
&\leq&\int_{S(r)}\log(|(f, H_{i})|^{2}+|(f,
H_{\sigma(i)})|^{2})^{\frac{1}{2}}\sigma\\&& +\int_{S(r)}\log(|(g,
H_{i})|^{2}+|(g, H_{\sigma(i)})|^{2})^{\frac{1}{2}}\sigma+O(1)\\
&\leq& T(r)+O(1).\end{eqnarray*} Therefore, for all $i\in\{1, 2,
\ldots, q\}$ we have
\begin{eqnarray*}&&T(r)+O(1)\\&\geq&N^{N}_{(f, H_{i}),
\leq m_{i}}(r)+N^{N}_{(f, H_{\sigma(i)}), \leq
m_{\sigma(i)}}(r)+N^{N}_{(g, H_{i}), \leq m_{i}}(r)+N^{N}_{(g,
H_{\sigma(i)}), \leq m_{\sigma(i)}}(r)\\&&-NN^{1}_{(f, H_{i}),\leq
m_{i}}(r)-NN^{1}_{(f, H_{\sigma(i)}),\leq
m_{\sigma(i)}}(r)+\sum_{j=1, j\neq i, \sigma(i)}^{q}N^{1}_{(f,
H_{j}), \leq m_{j}}(r).
\end{eqnarray*}
Note that $\sigma$ is bijective. Take summation of the above
inequality over $1\leq i\leq q,$ we have
\begin{eqnarray*}&&(q-2N-2)\sum_{j=1}^{q}N^{1}_{(f, H_{j}), \leq
m_{j}}(r)+2\sum_{j=1}^{q}\left(N^{N}_{(f, H_{j}), \leq
m_{j}}(r)+N^{N}_{(g, H_{j}), \leq m_{j}}(r)\right)\\&\leq&
qT(r)+O(1).\end{eqnarray*} By a similar discussion for $g$ instead
of $f,$ we have
\begin{eqnarray*}&&(q-2N-2)\sum_{j=1}^{q}N^{1}_{(g, H_{j}), \leq
m_{j}}(r)+2\sum_{j=1}^{q}\left(N^{N}_{(f, H_{j}), \leq
m_{j}}(r)+N^{N}_{(g, H_{j}), \leq m_{j}}(r)\right)\\&\leq&
qT(r)+O(1).\end{eqnarray*} Noting that $\frac{1}{N}N^{N}_{(f,
H_{j}), \leq m_{j}}(r)\leq N^{1}_{(f, H_{j}), \leq m_{j}}(r),$ we
get from the above inequalities that
\begin{eqnarray*}\frac{q+2N-2}{2N}\sum_{j=1}^{q}\left(N^{N}_{(f, H_{j}), \leq
m_{j}}(r)+N^{N}_{(g, H_{j}), \leq m_{j}}(r)\right)\leq q
T(r)+O(1).\end{eqnarray*}

By the Second Main Theorem, we get
\begin{eqnarray*}(q-N-1)T(r)&\leq&\sum_{i=1}^{q}\left(N^{N}_{(f, H_{i}), \leq
m_{i}}(r)+N^{N}_{(g, H_{i}), \leq
m_{i}}(r)\right)\\&&+\sum_{i=1}^{q}\left(N^{N}_{(f, H_{i}), \geq
m_{i}+1}(r)+N^{N}_{(g, H_{i}), \geq
m_{i}+1}(r)\right)+o(T(r)).\end{eqnarray*} Therefore, removing the
term $\sum_{i=1}^{q}\left(N^{N}_{(f, H_{i}), \leq
m_{i}}(r)+N^{N}_{(g, H_{i}), \leq m_{i}}(r)\right)$ from the above
inequalities we have
\begin{eqnarray*}&&\left(\frac{(q+2N-2)(q-N-1)}{2N}-q\right)T(r)\\&\leq&
\frac{q+2N-2}{2N}\sum_{j=1}^{q}\left(N^{N}_{(f, H_{j}), \geq
m_{j}+1}(r)+N^{N}_{(g, H_{j}), \geq
m_{j}+1}(r)\right)+o(T(r)).\end{eqnarray*} Noting that
\begin{eqnarray*}&&N^{N}_{(f, H_{j}), \geq m_{j}+1}(r)+N^{N}_{(g,
H_{j}), \geq m_{j}+1}(r)\\&\leq& \frac{N}{m_{j}+1}\left(N_{(f,
H_{j}), \geq m_{j}+1}(r)+N_{(g, H_{j}), \geq
m_{j}+1}(r)\right)\\&\leq& \frac{N}{m_{j}+1}\left(N_{(f,
H_{j})}(r)+N_{(g, H_{j})}(r)\right)\leq
\frac{N}{m_{j}+1}T(r)+o(T(r)),\end{eqnarray*} we get from the above
inequality that
\begin{eqnarray*}\left(\frac{(q+2N-2)(q-N-1)}{2N}-q\right)T(r)\leq
\frac{q+2N-2}{2}\sum_{j=1}^{q}\frac{1}{m_{j}+1}T(r)+o(T(r)).\end{eqnarray*}

Noting that $q=\sum_{j=1}^{q}\frac{m_{j}+1}{m_{j}+1},$ we can obtain
from the above inequality  that
\begin{eqnarray*}&&\left\|\left(\sum_{j=3}^{q}\frac{m_{j}}{m_{j}+1}-\frac{Nq-q+N+1}{N}+
\frac{4N-4}{q+2N-2}-\left(\frac{1}{m_{1}+1}+\frac{1}{m_{2}+1}\right)\right)T(r)\right.\\&&
\leq o(T(r)).\end{eqnarray*}
This is a contradiction.\vskip 6pt

\noindent{\bf Acknowledgements.} The authors would like to thank Dr.
Qi-Ming Yan, Dr. Hong-Zhe Cao and Dr. Kai Liu for making some
valuable suggestions to improve the present paper.\par

\vskip 6pt

\begin{thebibliography}{99}\bibitem{aihara}Y. Aihara, Unicity theorems for meromorphic
mappings with deficiencies, Complex Variables 42(2000), 259-268.
\bibitem{chen-yan-0} Z. H. Chen and Q. M. Yan, Uniqueness theorem of
meromorphic mappings into $\mathbb{P}^{N}(\mathbb{C})$ sharing
$2N+3$ hyperplanes regardless of multiplicities, Intern. J. Math.
20(2009), No. 6, 717-726.
\bibitem{dethloff-quang-tan} G. Dethloff, S. D. Quang and T. V. Tan,
A uniqueness theorem for meromorphic mappings with two families of
hyperplanes, Preprint.
\bibitem{dethloff-tan} G. Dethloff and T. V. Tan, Uniqueness
theorems for meromorphic mappings with few hyperplanes, Bull. Sci.
Math. 133(2009), 501-514.
\bibitem{dethloff-tan-1} G. Dethloff and T. V. Tan, An extension of uniqueness theorems for meromorphic
mappings, Vietnam J. Math. 34(2006), No. 1, 71-94.
\bibitem{fujimoto-1} H. Fujimoto, The uniqueness promble of meromorphic
maps into the complex projective spaces, Nagoya Math. J. 58(1975),
1-23.
\bibitem{fujimoto-2} H. Fujimoto, The uniqueness theorem for
algebraically non-degenerate meromorphic maps into
$\mathbb{P}^{N}(\mathbb{C}),$ Nagoya Math. J. 64(1976), 117-147.
\bibitem{hu-li-yang} P. C. Hu, P. Li and C. C. Yang, Unicity of
meromorphic mappings, Kluwer 2003.
\bibitem{nevanlinna}  R. Nevanlinna, Eindentig keitss\"{a}tze in der theorie
der meromorphen funktionen, Acta. Math.  48(1926), 367-391.
\bibitem{ru} M. Ru, Nevanlinna theory and its relation to
Diophantine Approximation, Singapore: World Scientific Publishing,
2001.
\bibitem{smiley} L. Smiley, Geometric conditions for unicity of
holomorphic curves, Contemp. Math. 25((1983), 149-154.
\bibitem{thai-quang} D. D. Thai and S. D. Quang, Uniqueness problem
with truncated multiplicities of meromorphic mappings in several
complex variables, Intern. J. Math. 17(2006), No. 10, 1223-1257.
\bibitem{yang} L. Yang, Multiple values of meromorphic functions and
functions combination, Acta Math. Sin. Chinese Ser. 14(1964),
428-437.(in Chinese)
\bibitem{yi-yang} H.-X. Yi and C.-C. Yang, Uniqueness theory of
meromorphic functions, Science Press 1995/Kluwer 2003.
\end{thebibliography}
\end{document}